\documentclass[graybox]{svmult}

\usepackage{mathptmx}       
\usepackage{helvet}         
\usepackage{courier}        
\usepackage{type1cm}        
\usepackage{makeidx}         
\usepackage{graphicx}        
\usepackage{multicol}        
\usepackage[bottom]{footmisc}

\usepackage{amsmath,amsfonts}

\makeindex             

\begin{document}

\title*{Optimality conditions and local regularity of the value function for the optimal exit time problem}

\author{Luong V. Nguyen}
\institute{Luong V. Nguyen \at Institute of Mathematics, Polish Academy of Sciences, Warsaw, Poland
\\ \email{vnguyen@impan.pl; luonghdu@gmail.com}
}
%
%
\maketitle

We consider the control problem with \textit{exit time}. Unlike the Bolza and Mayer problems, in this problem the terminal time of the trajectories is not fixed, but it is the first time at which they reach a given closed subset - \textit{the target}. The most studied example is the \textit{optimal time problem}, where we want to steer a point to the target in minimal time.

 In this section, we first introduce the exit time problem, then we recall the existence of optimal controls, and  some regularity results for the value function. We then use a suitable form of the \textit{Pontryagin maximum principle}  to study some optimality conditions and sensitivity relations for the exit time problem. The strongest regularity property for the value function that one can expect, in fairly general cases, is \textit{semiconcavity}\footnote{A function is semiconcave if it can be written as a sum of a concave function and a $C^2$ function.}. In this case, the value function is twice differentiable almost everywhere. Furthermore, in general, it fails to be differentiable at points where there are multiple  optimal trajectories and its differentiability at a point does not guarantee continuous differentiability around this point. In the subsection \ref{S3} we shown that, under suitable assumptions, the nonemptiness of proximal subdifferential of the value function at a point implies its continuous differentiability on a neighborhood of this point. 

\subsection{The optimal exit time problem} \label{S1}
We assume that a compact nonempty set $U\subset \mathbb{R}^m$ and a continuous function $f:\mathbb{R}^n \times U \to \mathbb{R}^n$ are given. We consider the \textit{control system}
\begin{equation}\label{Controlsys}
  \left\{
  \begin{array}{lcl}
    \dot{x}(t) & = & f(x(t),u(t)), \\
    x(0)& = & x_0 \in \mathbb{R}^n,
  \end{array}
  \right.\,\,\,\,\,\text{a.e.}\,\,t>0.
\end{equation}
where $u: \mathbb{R}_+ \to U$  is a measurable function which is called a \textit{control} for the system (\ref{Controlsys}).  The set $U$ is called the \textit{control set}. We denote by $\mathcal{U}_{ad}$ the set of all measurable control functions. We will often require the following assumptions
\begin{itemize}
\item[(A1)] \quad There exists $K_1>0$ such that 
$$|f(x_1,u) - f(x_2,u)| \le K_1|x_1-x_2|,\,\,\,\,\,\,\,\forall x_1,x_2 \in \mathbb{R}^n, u\in U.$$
\item[(A2)]\quad $D_xf$ exists and is continuous. Moreover, there exists $K_2 >0$ such that
$$||D_xf(x_1,u) - D_xf(x_2,u)|| \le K_2|x_1-x_2|, \,\,\,\,\,\,\,\forall x_1,x_2 \in \mathbb{R}^n, u\in U.$$
\end{itemize} 

It is well known from the ordinal differential equations theory that under assumption (A1), for each $u\in \mathcal{U}_{ad}$, (\ref{Controlsys}) has a unique solution. In this case, we will denote by $x^{x_0,u}(\cdot)$ the solution of (\ref{Controlsys}) and call $x^{x_0,u}(\cdot)$ the \textit{trajectory} (starting at $x_0$) of the control system (\ref{Controlsys}).

We now assume that a closed subset with compact boundary $\mathcal{K}$ of the state space $\mathbb{R}^n$ is given and is called the \textit{target}. For a given trajectory $x^{x_0,u}(\cdot)$ of (\ref{Controlsys}), we set
\begin{equation*}
\tau(x_0,u):=\min \left\{ t\ge 0: x^{x_0,u}(t) \in \mathcal{K}\right\},
\end{equation*}
with the convention that $\tau(x_0,u) = +\infty$ if $x^{x_0,u}(t)\not\in \mathcal{K}$ for all $t\ge 0$. Then $\tau(x_0,u)$ is the time at which the trajectory $x^{x_0,u}(\cdot)$ reaches the target for the first time, provided $\tau(x_0,u) < +\infty$ and we call $\tau(x_0,u)$ the \textit{exit time} of the trajectory $x^{x_0,u}(\cdot)$. Denote by $\mathcal{R}$ the set of all $x_0$ such that $\tau({x_0,u}) <+\infty$ for some $u(\cdot)\in \mathcal{U}_{ad}$ and we call $\mathcal{R}$ the \textit{reachable set}.

Given two continuous functions $L:\mathbb{R}^n \times U \to \mathbb{R}$ (called \textit{running cost}) and $\psi: \mathbb{R}^n \to \mathbb{R}$ (called \textit{terminal cost}) with $L$ positive and $\psi$ is bounded from below, we consider the functional
$$ J(x_0,u) = \int_0^{\tau(x_0,u)}{L\left(x^{x_0,u}(s),u(s) \right)}ds  + \psi\left( x^{x_0,u}(\tau(x_0,u))\right). $$
 We are interested in  minimizing $J(x_0,u)$, for $x_0\in \mathcal{R}$, over all $u(\cdot)\in \mathcal{U}_{ad}$. If $u^*(\cdot) \in \mathcal{U}_{ad}$ is such that 
$$J(x_0,u^*) = \min_{u\in \mathcal{U}_{ad}} J(x_0,u)$$ 
 then we call $u^*(\cdot)$ an \textit{optimal control} for $x_0$. In this case, $x^{x_0,u^*}(\cdot)$ is called an \textit{optimal trajectory}.

The \textit{value function} of the optimal exit time problem is defined by
\begin{equation*}
V(x_0): = \inf\left\{ J(x_0,u): u(\cdot) \in \mathcal{U}_{ad} \right\},\,\,\,\,\,\,\,\,\,x_0\in \mathcal{R}.
\end{equation*}
From the definition of $V$, we have the so-called \textit{dynamic programming principle}
\begin{equation*}
V(x_0) \le \int_0^t{L\left(x^{x_0,u}(s),u(s) \right)}ds + V\left( x^{x_0,u}(t)\right),\,\,\,\,\forall \,t\in [0,\tau(x_0,u)].
\end{equation*}
If $u(\cdot)$ is optimal then the equality holds.

The maximized Hamiltonian associated to the control system is defined by\footnote{Throughout the chapter, $p$ is usually  a row vector in $\mathbb{R}^{n,*}$. However, in this section, $\mathbb{R}^n$ and its dual space $\mathbb{R}^{n,*}$ are identical. We also denote by $a.b$ the inner product between two vectors $a$ and $b$.}
\begin{equation*}
\label{Ham}
\mathcal{H}(x,p) = \mathrm{max}_{u\in U} \left\{ -p.f(x,u) - L(x,u)\right\},\,\,\,\,\,\,\,\, (x,p)\in \mathbb{R}^n \times \mathbb{R}^n.
\end{equation*}
It is well-known, under some assumptions (see Theorem 8.18 in \cite{CSB04}), that $V$ is a viscosity solution of the Hamilton - Jacobi - Bellman equation
\begin{equation*}\label{HJB}
\mathcal{H}(x,\nabla V(x)) = 0.
\end{equation*}

We now list some more assumptions on the cost functionals and the target which will be used in the sequel.
\begin{itemize}
\item[(A0)]\quad  For all $x\in\mathbb{R}^n$, the following set is convex
$$\mathcal{F}(x):= \left\{(v,\lambda)\in\mathbb{R}^{n+1}: \exists u\in U\,\,\text{such that}\,\, v= f(x,u),\lambda \ge L(x,u)   \right\}  .$$
\item[(A3)]\quad There exist $N>0$ and $\alpha>0$ such that $|f(x,u)| \le N$ and $L(x,u) \ge \alpha$ for all $x\in \mathbb{R}^n$ and $u\in U$.
\item[(A4)]\quad The function $L$ is continuous in both arguments and locally Lipschitz continuous with respect to $x$, uniformly in $u$. Moreover, $L_x(x,u)$ exists for all $x,u$ and is locally Lipschitz continuous in $x$, uniformly in $u$.
\item[(A5)]\quad There exists a neighborhood $\mathcal{N}$ of $\text{bdry}\mathcal{K}$ such that $\psi$ is locally semiconcave and is of class $C^1$ in $\mathcal{N}$. Moreover, denoting by  $G$ the Lipschitz constant of $\psi$ in $\mathcal{N}$, we assume 
$$
G< \frac{\alpha}{N}.
$$
\item[(A6)]\quad The boundary of $\mathcal{K}$ is an $(n-1)$-dimensional manifold of class $C^{1,1}_{loc}$ and there exists $\gamma >0$ such that for any $z \in \text{bdry}\mathcal{K}$, we have
$$ \min_{u\in U} \,\, f(z,u).n_z \le -\gamma,$$
where $n_z$ denotes the unit outward normal to $\mathcal{K}$ at $z$.
\end{itemize}
Assumption (A0) is a condition to ensure the existence of optimal trajectories. More precisely, one has
\begin{theorem} \label{Existence} \cite{CPS00,CSB04}
Under assumptions (A0) - (A5), there exists a minimizer for optimal control problem for any choice of initial point $y\in \mathcal{R}$. Moreover, the uniform limit of optimal trajectories is an optimal trajectory; that is, if $x_k(\cdot)$ are trajectories converging uniformly to $x(\cdot)$ and every $x_k(\cdot)$ is optimal for the point $y_k: = x_k(0)$, then $x(\cdot)$ is optimal for $y: = \lim y_k$.
\end{theorem}
The condition $G<\alpha /N$ in assumption (A5) can be regarded as a compatibility condition on the terminal cost $\psi$. Together with other assumptions, it ensures the continuity of the value function (see Remark 2.6 in \cite{CPS00} and Proposition IV.3.7 in \cite{BCD97}). Furthermore, we have the following regularity property of the value function.
\begin{theorem} \cite{CPS00,CSB04}
Under hypothesis (A1)-(A6), the value function $V$ is locally semiconcave in $\mathcal{R}\setminus \mathcal{K}$.
\end{theorem}
Note that in \cite{CPS00,CSB04}, the semiconcavity result is proved under weaker assumptions on the data. In fact, $\mathcal{K}$ is only assumed to satisfy an interior sphere condidtion, while $f$, $L$ and $\psi$ are assumed to be semiconcave in the $x$-variable and $L_x$ is only continuous. \\
 For  the precise definition, properties and characterizations of semiconcave functions, we refer the reader to \cite{CSB04}.
\subsection{Optimality conditions and sensitivity relations}\label{S2}
We present some optimality conditions and sensitivity relations for the optimal exit time problem. One of important tools for our analysis is given by the so-called \textit{Pontryagin maximum principle}. Before recalling a version of the maximum principle for the optimal control problem under consideration, we need to introduce some notation. For a given subset $A$ of $\mathbb{R}^n$, we denote by $\text{bdry}A$ its boundary, by $A^c$ its complement. The distance function from $A$ is define for $x\in \mathbb{R}^n$ as
$$d_A(x) := \inf_{y\in A}|x-y|,\quad x\in \mathbb{R}^n,$$
and the oriented distance function from $A$ is defined by $b_A(x) := d_A(x) - d_{A^c}(x), x\in \mathbb{R}^n$, whenever $A\ne \mathbb{R}^n$.

Let $\Omega$ be an open subset of $\mathbb{R}^n$, $h:\Omega \to \mathbb{R}$ be a lower semicontinuous function and $x\in \Omega$. The \textit{proximal subdifferential} of $h$ at $x$ is the set
\begin{eqnarray*}
\partial^Ph(x) &:=& \left\{v\in \mathbb{R}^n: \text{there exist}\,\, c>0,\rho >0\,\,\text{such that}\right. \\
&&\left. \qquad h(y) - h(x) - v.(y-x) \ge - c|y-x|^2,\quad \forall y\in B(x,\rho)\right\}.
\end{eqnarray*}
The \textit{Fr\'echet subdifferential} of $h$ at $x$ is the set 
$$ D^-h(x): = \left\{v\in \mathbb{R}^n: \liminf_{y\to x} \frac{h(y) - h(x) -   v.(y-x)}{|y-x|}    \ge 0 \right\}.   $$
The \textit{Fr\'echet superdifferential} of $h$ at $x$ is the set 
$$ D^+h(x): = \left\{v\in \mathbb{R}^n: \limsup_{y\to x} \frac{h(y) - h(x) -  v.(y-x)}{|y-x|}    \le 0 \right\}.   $$
If $h$ is locally Lipschitz, then the reachable gradient of $h$ at $x$ is the set
$$D^* h(x) := \left\{v\in \mathbb{R}^n: \exists x_n \to x, \nabla h(x_n) \to v\,\,\text{as}\,\,n\to\infty     \right\}.$$

We now start with two technical lemmas.
\begin{lemma} (see, e.g. \cite{CPS00}) \label{TechLem1}
Assume (A1) - (A6). Given $z\in \text{bdry}\mathcal{K}$, let $\zeta$ be the outer normal to $\mathcal{K}$ at $z$. Then there exists a unique $\mu>0$ such that $H(z,\nabla \psi(z) +\mu\zeta) = 0$.
\end{lemma}
Notice that since the boundary of the target $\mathcal{K}$ is of class $C^{1,1}_{loc}$, the outer normal to $\mathcal{K}$ at a point $z\in \text{bdry}\mathcal{K}$ is $\nabla b_{\mathcal{K}}(z)$.  From Lemma \ref{TechLem1}, the function $\mu: \text{bdry}\mathcal{K} \to \mathbb{R}_{+}$ which satisfies $H(z, \nabla \psi(z) + \mu(z) \nabla b_{\mathcal{K}}(z)) = 0$ is well-defined. Moreover, we have
\begin{lemma}\label{TechLem2} \cite{LN}
Assume (A1) - (A6). The function $\mu: \text{bdry}\mathcal{K} \to \mathbb{R}_+$ is continuous.
\end{lemma}

We recall the maximum principle in the following form

\begin{theorem} \label{MPM}
Assume (A1) - (A6). Let $x\in \mathcal{R}\setminus \mathcal{K}$ and let $\bar{u}$ be an optimal control for $x_0$. Set for simplicity
$$x(t): = x^{x_0,\bar{u}}(t),\,\,\,\,\,\tau:= \tau(x_0,\bar{u}),\,\,\,\,\,\,\,z:= x({\tau}).$$
Let $p\in W^{1,1}([0,\tau]; \mathbb{R}^n)$ be the solution to the equation
\begin{equation}\label{ADE}
\dot{p}(t)=D_xf(x(t),\bar{u}(t))^\top p(t) - L_x(x(t),\bar{u}(t)),
\end{equation}
with $
    p(\tau) =  \nabla \psi(z) + \mu(z)\nabla b_{\mathcal{K}}(z).$\\
Then $p$ satisfies
$$ - p(t).f(x(t),\bar{u}(t))- L(x(t),\bar{u}(t)) = \mathcal{H}(x(t),p(t)),  $$
for a.e. $t\in [0,\tau]$
\end{theorem}
For the proof of the above maximum principle, we refer the reader to Theorem 4.3 in \cite{CPS00} where the principle is proved under weaker assumptions on $L$ and $\psi$.

Given an optimal trajectory $x(\cdot)$, then, by Lemma \ref{TechLem1},  there is a unique function $p(\cdot)$ satisfying the properties of Theorem \ref{MPM} and we call $p(\cdot)$ the \textit{dual arc} associated to the trajectory $x(\cdot)$. Observe that the dual arc is a nonzero function and satisfies $p(\tau) = \nabla \psi(x(\tau)) + \mu(x(\tau))\nabla b_{\mathcal{K}}(x(\tau))$ where $\tau$ is the exit time of $x(\cdot)$. 
The following theorem gives a connection between the dual arcs and the \textit{Fr\'echet supperdifferential} of the value function.
\begin{theorem}\label{DuSup} \cite{CPS00}
Under the assumptions of Theorem \ref{MPM}, the dual arc $p(\cdot)$ satisfies
$$p(t) \in D^+V(x(t)),\,\,\,\,\,\,\,\,\forall \,t\in [0,\tau).$$
\end{theorem}
It is proved in \cite{CPS00,CSB04} that under the assumptions of Theorem \ref{MPM} and the following assumption
\begin{itemize}
\item[(H)] \quad for any $x\in \mathbb{R}^n$, if $\mathcal{H}(x,p) = 0$ for all $p$ in a convex set $C$, then $C$ is a singleton,
\end{itemize}
 the value function $V$ is differentiable along optimal trajectories except the initial and final point points and therefore, by Theorem \ref{DuSup} if $p(\cdot)$ is the dual arc associated  with an optimal trajectory $x^{x_0,u}(\cdot)$ then $p(t) = \nabla V(x^{x_0,u}(t))$ for all $t\in (0,\tau(x_0,u))$. This property plays an important role to prove an one-to-one correspondence between the number of optimal trajectories starting at a point $x_0\in \mathcal{R} \setminus \mathcal{K}$ and the number of elements of the \textit{reachable gradient} $D^*V(x_0)$ of $V$ at $x_0$. This implies that $V$ is differentiable at $x_0$ iff there is a unique optimal trajectory starting at $x$.  The following example shows that without assumption (H), $V$ may not differentiable along optimal trajectories.
 \begin{example} \label{EX1}
We consider the minimum time problem i.e., $L\equiv 1, g\equiv 0$, for the control system
$$
\left\{
  \begin{array}{lcl}
    \dot{x}_1(t) &=& u\\
     \dot{x}_2(t) &=& 0
  \end{array},\,\,\,\,u\in U:=[-1,1],
  \right.
$$
with the initial conditions $x_1(0) = y_1, x_2(0) = y_2$. Define the set
\begin{eqnarray*}
\mathcal{D} &=& \{(y_1,y_2)^\top: 2y_1-3y_2-2 >0\} \cap \{(y_1,y_2)^\top: 2y_1+3y_2-2 >0\}\\
&\cap& \{(y_1,y_2)^\top: 2y_1+3y_2-14 <0\} \cap  \{(y_1,y_2)^\top: 2y_1-3y_2-14 <0\}.
\end{eqnarray*}
The target is the set
$$\mathcal{K} = \mathbb{R}^2 \setminus \mathcal{D}.$$
The Hamiltonian is
$$\mathcal{H}(x,p) = \sup_{u\in U} \left\{- \begin{pmatrix}
u\\0
\end{pmatrix}. \begin{pmatrix}
p_1\\p_2
\end{pmatrix}  -1 \right\} = |p_1| - 1,\,\,\,\forall x\in \mathbb{R}^2, p=(p_1,p_2)^\top\in \mathbb{R}^2.
$$
One can easily check that all assumptions of Theorem 3.8 in \cite{CPS00} (which says that the value function is differentiable along optimal trajectories except the starting and the terminal points) are satisfied and assumption (H) is not satisfied. Let $T(\cdot)$ be the minimum time to reach the target.\\
If $y = (y_1,y_2)^\top\in \mathcal{D} \cap \{(y_1,y_2)^\top\in \mathbb{R}^2: y_1 < 4\}$, then $u^*(\cdot) \equiv -1$ is the optimal control for $y$ and 
$$x(t) = \begin{pmatrix}
y_1 - t\\
y_2
\end{pmatrix},\,\,\,\text{for all}\,\,t\in [0,T(y)]$$
is the optimal trajectory starting at $y$ and we can easily compute that
$$T(y) = y_1 - \frac{3}{2}|y_2|-1.$$
If $y = (y_1,y_2)^\top\in \mathcal{D} \cap \{(y_1,y_2)^\top\in \mathbb{R}^2: y_1 > 4\}$, then $u^*(\cdot) \equiv 1$ is the optimal control for $y$, the optimal trajectory is
$$x(t) = \begin{pmatrix}
y_1 + t\\
y_2
\end{pmatrix},\,\,\,\text{for all}\,\,t\in [0,T(y)],$$
and the minimum time to reach the target from $y$ is 
$$T(y) = -y_1 - \frac{3}{2}|y_2| + 7.$$
Since $T$ is not differentiable when $y_2 = 0$, $T$ fails to be differentiable at any point of optimal trajectories starting at $(y_1,0)^\top \in \mathcal{D}$.
\end{example}

Later we will see that $V$ is still differentiable at a point $x$ iff there is a unique optimal trajectory starting at $x$ even when (H) is not satisfied. In this case we may not have an one-to-one correspondence between the number of optimal trajectories starting at a point $x\in \mathcal{R} \setminus \mathcal{K}$ and the number of elements of the reachable gradient $D^*V(x)$.

If we assume that
\begin{itemize}
\item[(H1)]\quad $\mathcal{H} \in C^{1,1}_{loc}\left(\mathbb{R}^n \times (\mathbb{R}^n\setminus \{0\})\right)$
\end{itemize}
then we can compute partial derivatives of the maximized Hamiltonian (see Theorem 7.3.6 and also Remark 8.4.11 in \cite{CSB04}).
\begin{theorem} \label{DeHam}
If (H1) holds, then for any $(x,p) \in \mathbb{R}^n\times (\mathbb{R}^n\setminus \{0\})$, we have
$$\mathcal{H}_p(x,p) = - f(x,u^*(x,p)),$$
and
$$\mathcal{H}_x(x,p) = -D_xf(x,u^*(x,p))^\top p - L_x(x,u^*(x,p)),$$
where $u^*(x,p)$ is any element of $U$ such that
$$ -f(x,u^*(x,p)).p - L(x,u^*(x,p)) = \mathcal{H}(x,p).$$
\end{theorem}
Since we are going to evaluate the Hamiltonian $\mathcal{H}$ along dual arcs which are nonzero, the lack of differentiability of $\mathcal{H}$ is not an obstacle. From Theorem \ref{MPM} and Theorem \ref{DeHam}, we have
\begin{theorem}\label{HamSys}
Assume (A1) - (A6) and (H1). Let $x(\cdot)$ be an optimal trajectory and let $p(\cdot)$ be the associated dual arc to $x(\cdot)$. Then the pair $(x(\cdot),p(\cdot))$ solves the system
\begin{equation}\label{Hsys}
\left\{
  \begin{array}{lcl}
\dot{x}(t)&=&-\mathcal{H}_p(x(t),p(t)) \\
 \dot{p}(t)& = & \mathcal{H}_x(x(t),p(t)).
  \end{array}
  \right.
\end{equation}
Consequently $x(\cdot)$ and $p(\cdot)$ are of class $C^1$.
\end{theorem}
The next theorem can be seen as a propagation property of the \textit{Fr\'echet subdifferential} of the value function forward in time along optimal trajectories
\begin{theorem}\label{PropaSub}
Assume (A1), (A2) and (A4). Let $x_0 \in \mathcal{R}\setminus \mathcal{K}$ and  let $\bar{u}(\cdot)$ be an optimal control for $x_0$. Set for simplicity
$$\bar{x}(t):=x^{x_0,\bar{u}}(t),\,\,\,\,\,\,\,\,\,\,\tau:= \tau(x_0,\bar{u}). $$
Assume that $D^- V(x_0)\ne \emptyset$ and let $p \in W^{1,1}( [0,\tau];\mathbb{R}^n)$ be a solution of the equation (\ref{ADE})
satisfying $p(0) \in D^-V(x_0)$. Then $p(t) \in D^-V(\bar{x}(t))$ for all $t\in [0,\tau)$.
\end{theorem}
For the proof of the previous Theorem, one can find in \cite{LN}. Similarly, one can prove the following propagation result for the proximal subdifferential of the value function which will be used to prove the main results in the next section.
\begin{theorem}\label{PropaProSub}
Assume (A1), (A2) and (A4). Let $x_0 \in \mathcal{R}\setminus \mathcal{K}$ and  let $\bar{u}(\cdot)$ be an optimal control for $x_0$. Set for simplicity
$$\bar{x}(t):=x^{x_0,\bar{u}}(t),\,\,\,\,\,\,\,\,\,\,\tau:= \tau(x_0,\bar{u}). $$
Assume that $\partial^P V(x_0)\ne \emptyset$ and let $p\in W^{1,1}([0,\tau];\mathbb{R}^n)$ be a solution of the equation (\ref{ADE})
satisfying $p(0) \in \partial^PV(x_0)$. Then for some $c>0$ and for all $t\in [0,\tau)$, there exists $r>0$ such that, for every $z\in B(\bar{x}(t),r)$,
$$V(z) - V(\bar{x}(t)) \ge  p(t).(z - \bar{x}(t)) - c|z-\bar{z}(t)|^2.$$
Consequently, $p(t) \in \partial^PV(\bar{x}(t))$ for all $t\in [0,\tau)$.
\end{theorem}
Using above results, we can obtain the following results 9see \cite{LN} for the proofs).
\begin{theorem}\label{UniqTr}
Assume (A0) - (A6) and (H1). Let $x_0 \in \mathcal{R} \setminus \mathcal{K}$ be such that $V$ is differentiable at $x_0$. Consider the solution $(x(\cdot),p(\cdot))$ of (\ref{Hsys}) with initial conditions
\begin{equation*}\label{Ini1}
\left\{
  \begin{array}{lcl}
x(0))&=&x_0 \\
 p(0)& = & DV(x_0).
  \end{array}
  \right.
\end{equation*}
Then $x(\cdot)$ is an optimal trajectory for $x_0$ and $p(\cdot)$ is the dual arc associated to $x(\cdot)$ with $p(t) = DV(x(t))$ for all $t\in [0,\tau)$ where $\tau$ is the exit time of $x(\cdot)$. Moreover, $x(\cdot)$ is the unique optimal trajectory staring at $x_0$.
\end{theorem}
\begin{theorem} \label{ReGra}
Assume (A0) - (A6) and (H1). Let $x_0\in \mathcal{R}\setminus \mathcal{K}$ and $q\in D^*V(x_0)$. Consider the solution $(x(\cdot),p(\cdot))$ of (\ref{Hsys}) with initial conditions
\begin{equation}\label{Ini2}
\left\{
  \begin{array}{lcl}
x(0))&=&x_0 \\
 p(0)& = & q.
  \end{array}
  \right.
\end{equation}
Then $x(\cdot)$ is an optimal trajectory for $x_0$ and $p(\cdot)$ is the dual arc associated to $x(\cdot)$. Moreover $p(t) \in D^*V(x(t))$ for all $t\in [0,\tau)$ where $\tau$ is the exit time of $x(\cdot)$.
\end{theorem}
\begin{theorem} \label{UniOpt}
Assume (A0) - (A6) and (H1). If there is only one optimal trajectory starting at a  point $x\in \mathcal{R}\setminus \mathcal{K}$ then $V$ is differentiable at $x$.
\end{theorem}
From Theorem \ref{UniqTr} and Theorem \ref{UniOpt}, we have
\begin{corollary} \label{Coro2}
Assume (A0) - (A6) and (H1). The value function $V$ is differentiable at a point $x\in \mathcal{R}\setminus \mathcal{K}$ if and only if there exists a unique optimal trajectory starting at $x$.
\end{corollary}
In Example \ref{EX1}, the value function is not differentiable at any point $(y_1,0)^\top \in \mathcal{D}$ although there is a unique optimal trajectory for every $(y_1,0)^\top$ with $y_1 \ne 4$. The reasons are that the maximized Hamiltonian $\mathcal{H}$ does not belong to $C^{1,1}_{loc}(\mathbb{R}^2 \times (\mathbb{R}^2\setminus \{0\}))$ and that the target is not smooth. We now give a simple example showing that the value function is differentiable at a point $x$ although there are multiple optimal trajectories starting at $x$.
\begin{example} \label{EX2} We consider the minimum time problem for the control system
\begin{equation*}
\begin{pmatrix}
\dot{y}_1(t)\\
\dot{y}_2(t)
\end{pmatrix}
= \begin{pmatrix}
u_1\\
u_2
\end{pmatrix},\,\,\,\,\,\,\,|u_i|\le 1,\,\,i=1,2,
\end{equation*}
with the initial condition $y_1(0) = x_1, y_2(0) = x_2$. The target is the set
\begin{eqnarray*}
 \mathcal{K} &=&\left\{(x_1,x_2)^\top : x_1 \le 0 \right\} \cap \left\{(x_1,x_2)^\top: x_2\le 4 + \sqrt{-x_1^2 - 4x_1} \right\}\\
 &\cap& \left\{(x_1,x_2)^\top: x_1\ge -4 \right\} \cap \left\{(x_1,x_2)^\top: x_2\ge -4 - \sqrt{-x_1^2 - 4x_1} \right\}
\end{eqnarray*}
The Hamiltonian is defined by $\mathcal{H}(x,p) = |p_1| + |p_2| - 1, \forall x \in \mathbb{R}^2, p=(p_1,p_2)^\top \in \mathbb{R}^2$. We can easily check that assumptions (A0) - (A6) are satisfied anh (H1) is not satisfied.

Observe that the minimum time function (the value function) $T$ is of class $C^{1,1}_{loc}(\mathcal{R}\setminus \mathcal{K})$ (see, e.g.,  \cite{CS95,CSB04}). Therefore $T$ is differentiable at $x=(1,0)^\top$. However, there are multiple optimal trajectories starting at $x$. Indeed, the trajectories corresponding to the controls $u_1(\cdot) \equiv (-1,0)^\top, u_2(\cdot)  \equiv (-1,1)^\top$ and $u_3(\cdot) \equiv (-1,-1)^\top$ are optimal for $x$.
\end{example}
We next give a class of control systems which can be applied Corollary \ref{Coro2}.
\begin{example} [see, e.g. Example 4.12 \cite{CPS00}] \label{EX3}
We consider the control system with the dynamics given by
$$f(x,u) = h(x) +\sigma (x)u,$$
where $h:\mathbb{R}^n \to \mathbb{R}^n$, $\sigma: \mathbb{R}^n \to \mathbb{M}^{n\times n}$ and the control set $U$ is the closed ball of center zero and radius $R>0$ in $\mathbb{R}^n$. We also consider the running cost of the form
$$L(x,u) = \ell (x) + \frac{1}{2}|u|^2,$$
where $\ell: \mathbb{R}^n \to \mathbb{R}$.

Since $f$ is affine and $L$ is convex with respect to $u$ and $U$ is convex, one can check that assumption (A0) is satisfied.  If we assume that $\sigma, h, \ell$ are  of class $C^{1,1}$, that $\sigma, h$ are bounded and Lipschitz and that $\ell$ is bounded below by a positive constant, then assumption (A1) - (A4) are satisfied. The Hamiltonian
\begin{eqnarray*}
\mathcal{H}(x,p) &=& \max_{u\in U} \left\{  - (h(x) + \sigma(x)u).p - l(x) - \frac{1}{2}|u|^2   \right\}\\
  &=&\left\{
  \begin{array}{lcl}
     -h(x).p  -\ell(x) + \frac{1}{2}|\sigma(x)^\top p|^2 &\,\,\,\,\, \text{if}\,\, |\sigma(x)^\top p| \le R\\
     -h(x).p -\ell(x) + R |\sigma(x)^\top p|- \frac{R^2}{2} &\,\,\,\,\, \text{if}\,\, |\sigma(x)^\top p| > R
  \end{array}
  \right.
\end{eqnarray*}
satisfies assumption (H1). Then if the final cost function $\psi$ and the target satisfy assumption (A5) and  (A6) then our result can be applied. 
\end{example}
\subsection{Local regularity of the value function}\label{S3}
In this section, we provide sufficient conditions which guarantee the continuous differentiability of the value function $V$ around a given point.  Local $C^1$ regularity of $V$ is discussed in the subsection \ref{S1}, whereas local $C^k$ ($k\ge 2$) regularity of $V$ is established in the subsection \ref{S2}. In both subsections \ref{S1} and \ref{S2}, the main condition to ensure the continuous differentiability of $V$ around a given point $x$ is the nonemptiness of the proximal subdifferential of $V$ at $x$.
 
\subsubsection{Local $C^1$ regularity}\label{S1}
In addition, we require the following assumptions.
\begin{itemize}
\item[(A7)]\quad $\psi$ is of class $C^2$ in a neighborhood of bdry$\mathcal{K}$ and bdry$\mathcal{K}$ is of class $C^2$ .
\item[(H2)]\quad $\mathcal{H} \in C^2_{loc}\left(\mathbb{R}^n \times (\mathbb{R}^n\setminus \{0\})\right)$.
\end{itemize} 

Below we denote by  $T_{\text{bdry}\, \mathcal{K}}(z)$ the tangent
space to the $(n-1)-$dimensional $C^2$-manifold  $\text{bdry}\,
\mathcal{K}$ at $z \in \text{bdry}\, \mathcal{K}$.

Consider the \textit{Hamiltonian system}
\begin{equation}\label{Hamsys}
  \left\{
  \begin{array}{lcl}
    -\dot{x}(t) & = & \mathcal{H}_p(x(t),p(t)) \\
    \,\,\dot{p}'(t)& = & \mathcal{H}_x(x(t),p(t)),
  \end{array}
  \right.
\end{equation}
on $[0,T]$ for some $T > 0$, with the final conditions
\begin{equation}\label{FinCd}
  \left\{
  \begin{array}{lcl}
    x(T) & = & z \\
    p(T)& = & \varphi(z),
  \end{array}
  \right.
\end{equation}
where $z$ is in a neighborhood of $\text{bdry}\,\mathcal{K}$ and
$\varphi(z) = \nabla \psi(z)+ \mu(z) \nabla b_{\mathcal{K}}(z)$ with $\mu(\cdot)$
satisfying $\mathcal{H}(z,\nabla \psi(z) + \mu(z)\nabla b_{\mathcal{K}}(z)) = 0$. Note that, by (A7),  $\mu(\cdot)$ is of class $C^1$ in a neighborhood bdry$\mathcal{K}$ (see Proposition 3.2 in \cite{P02}) and therefore  $\varphi(\cdot)$ is of class $C^1$ in a neighborhood of bdry$\mathcal{K}$.

For a given $z$ in a neighborhood of
$\text{bdry}\,\mathcal{K}$, let $(x(\cdot;z),p(\cdot;z))$
be the solution of (\ref{Hamsys}) - (\ref{FinCd}) defined
on a time interval $[0,T]$ with  $T>0$. Consider the
so-called \textit{variational system}
\begin{equation}\label{Varsys}
  \left\{
  \begin{array}{lcl}
    -\dot{X} & = & \mathcal{H}_{xp}(x(t),p(t))X + \mathcal{H}_{pp}(x(t),p(t))P,\,\,\,\,\,\,X(T) = I \\
    \,\,\,\dot{P}& = & \mathcal{H}_{xx}(x(t),p(t)) X + \mathcal{H}_{px}(x(t),p(t))P,\,\,\,\,\,\,P(T) =
    D\varphi(z).
  \end{array}
  \right.
\end{equation}
Then the solution $(X,P)$ of (\ref{Varsys}) is defined in
$[0,T]$ and depends on $z$. Moreover
$$ X(\cdot;z) = D_zx(\cdot;z)\,\,\,\,\,\text{and}\,\,\,\, P(\cdot;z) = D_zp(\cdot;z), \qquad \text{on}\,\,[0,T].$$
\begin{definition}\label{Conj}
For $z\in \text{bdry}\, \mathcal{K}$, the time
$$ t_c(z): = \inf\{t\in [0,T]: X(s)\theta \ne 0, \forall\;  0 \neq\theta \in T_{\text{bdry}\, \mathcal{K}}(z), \forall s \in [t,T] \}$$
is said to be \textit{conjugate-like} for $z$ iff there exists $0 \neq\theta \in
T_{\text{bdry}\, \mathcal{K}}(z)$ such that
$$X(t_c(z))\theta = 0. $$
In this case, the point $x(t_c(z))$ is called conjugate-like for $z$.
\end{definition}
\begin{remark} In the classical definition of conjugate point it is required, for some $0 \neq \theta \in \mathbb{R}^N$,  $X(t_c(z)) \theta = 0$
 (see e.g. \cite{CS951,CaF,P02} and Section \ref{S2}).  Here, we have narrowed the
set of such  $\theta$ getting  then  a stronger result in Theorem
\ref{main1} below than the one we would have with the classical
definition.
\end{remark}

\begin{theorem}\label{main1}
Assume (A0) - (A7) and (H2). Let $x_0 \in \mathcal{R}\setminus
\mathcal{K}$ be such that $V$ is differentiable at
$x_0$ and $x^{x_0,u}(\cdot)$ be the optimal trajectory for $x_0$.
Set $z = x^{x_0,u}(\tau(x_0,u))$.  If there is no conjugate-like
time in $[0,\tau(x_0,u)]$ for $z$ then
$V$ is of class $C^1$ in a neighborhood of $x_0$.
\end{theorem}

When the maximized Hamiltonian is strictly convex with respect to the second variable, we can progress as in \cite{CF1,CF2} to obtain the following result.
\begin{theorem} \label{main2}
Assume  (A0) - (A7),  (H2)  and that $\mathcal{H}_{pp}(x,p)>0$ for all $(x,p) \in \mathbb{R}^n \times
(\mathbb{R}^n \setminus \{0\})$. Let $x_0 \in \mathcal{R} \setminus
\mathcal{K}$. If $\partial^P V(x_0) \ne \emptyset$,
then $V$ is of class $C^1$ in a neighborhood of
$x_0$.
\end{theorem}

Since $V$ is locally semiconcave and $\partial^P V(x_0) \ne \emptyset$ , $V$ is differentiable at $x_0$. The idea of the proof is to absent a conjugate time for the final point of the optimal trajectory starting at $x_0$ and then apply Theorem \ref{main1}.

In Example \ref{EX3}, if $\sigma(x)$ is nonsingular for all $x\in \mathbb{R}^n$ then $\mathcal{H}_{pp}(x,p) >0$ for all $(x,p) \in \mathbb{R}^n \times
(\mathbb{R}^n \setminus \{0\})$. Therefore if $h,\sigma,\ell,\psi$ and bdry$\mathcal{K}$ are smooth enough then Theorem \ref{main2} can be applied.  

When the running cost does not depend on $u$, i.e., $L = L(x)$, the maximized Hamiltonian is never strictly convex with respect to the second variable. In this case $0\ne p \in \ker \mathcal{H}_{pp}(x,p)$ for all $x\in \mathbb{R}^n$ whenever $\mathcal{H}_{pp}(x,p)$ exists. Following the lines for the minimum time problem in \cite{FL13}, we obtain the following particular case
\begin{theorem} \label{main3}
Assume that  (A0) - (A7),  (H2) hold true,  the kernel of $\mathcal{H}_{pp}(x,p)$
has the dimension equal to $1$ for every $(x,p) \in \mathbb{R}^n \times
(\mathbb{R}^n \setminus \{0\})$ and that $L = L(x)$. Let $x_0 \in \mathcal{R} \setminus
\mathcal{K}$. If $\partial^P V(x_0) \ne \emptyset$,
then $V$ is of class $C^1$ in a neighborhood of
$x_0$.
\end{theorem}

\begin{example} [see, e.g., \cite{FL13}] \label{EX4}
Consider the control system with the dynamics given by
$$f(x,u) = h(x) + \sigma(x)u,$$
where $h:\mathbb{R}^n \to \mathbb{R}^n, \;\sigma:\mathbb{R}^n \to \mathcal{L}(\mathbb{R}^n; \mathbb{R}^n)$ and
the control set $U$ is the closed ball in $\mathbb{R}^n$ of center zero
and radius
$R>0$.\\
Since $f$ is affine with respect to $u$, assumption (A0) is
verified. Let $L(x,u) = L(x)$ for all $(x,u)\in \mathbb{R}^n\times U$. The Hamiltonian
\begin{eqnarray*}
\mathcal{H}(x,p) &=&\max_{u\in U} \left\{  -(h(x) + \sigma(x)u).p   \right\} - L(x)\\
&=& - h(x).p + \max_{u\in U} \left\{ -u.\sigma(x)^\top p\right\}- L(x)\\
&=& - h(x).p  + |\sigma(x)^\top p| - L(x)
\end{eqnarray*}
satisfies assumption (H2) whenever $\sigma(x)$ is also surjective for
all $x\in \mathbb{R}^n$ and $h,\sigma, L$ are of class $C^2$.  Furthermore, for all $(x,p)\in \mathbb{R}^n
\times (\mathbb{R}^n\setminus \{0\})$
$$\mathcal{H}_p(x,p) = -h(x) + \frac{1}{|\sigma(x)^\top p|} \sigma(x)\sigma(x)^\top p$$
and  for any $q \in \mathbb{R}^n$,
$$\mathcal{H}_{pp}(x,p)(q,q) = \frac{1}{|\sigma(x)^\top p|}\; |\sigma(x)^\top q|^2 - \frac{1}{|\sigma(x)^\top p|^3}\left( \sigma(x)^\top p. \sigma(x)^\top q \right) ^2.$$
Fix any  $q \in $ker $\mathcal{H}_{pp}(x,p)$. Then, from the above equality
we get
\begin{equation}\label{march29}
|\sigma(x)^\top p|^2 |\sigma(x)^\top q|^2 =  \left(\sigma(x)^\top p. \sigma(x)^\top q  \right)^2 .
\end{equation}
On the other hand, if $\sigma(x)^\top q \notin \mathbb{R} \left(\sigma(x)^\top p \right)$,
then $$| \sigma(x)^\top p. \sigma(x)^\top q | < |\sigma(x)^\top p| |\sigma(x)^\top q|.$$
Hence, by (\ref{march29}),  $\sigma(x)^\top q \in \mathbb{R} \sigma(x)^\top p$. Let $\lambda
\in \mathbb{R}$ be such that $\sigma(x)^\top q = \lambda \sigma(x)^\top p$. Consequently
$\sigma(x)^\top (q- \lambda p)=0$.  Since $\sigma(x)$ is surjective, we deduce
that $q= \lambda p$ and that $ q \in  \mathbb{R} p$.

Using the inclusion $p \in $ker $\mathcal{H}_{pp}(x,p)$, we deduce that
$\ker \mathcal{H}_{pp}(x,p) = \mathbb{R} p$ for all $(x,p) \in \mathbb{R}^n\times
(\mathbb{R}^n\setminus \{0\})$, i.e., $$\dim \ker \mathcal{H}_{pp}(x,p) = 1, \;\;
\forall \;  (x,p) \in \mathbb{R}^n\times (\mathbb{R}^n\setminus \{0\}).$$ So, if
the target $\mathcal{K}$ and $\psi$ are of class  $C^2$ and for any $z \in
$bdry $\mathcal{K}$, the classical inward pointing condition
$$\min_{u \in U} \,\, (h(z) + \sigma(z)u).n_z <0$$
holds true, then  Theorem \ref{main3} can be applied. 
\end{example}
\subsubsection{Local $C^k$  regularity} \label{S2}
Let $k$ be an integer with $k\ge 2$. In this subsection, we require the following additional assumptions.
\begin{itemize}
\item[(A8)]\quad The functions $f$ and $L$ are of class $C^k$ in both arguments and the boundary of the target $\mathcal{K}$ is an $(n-1)$- manifold of class $C^{k+1}$. Moreover, $\psi$ is of class $C^{k+1}$ in a neighborhood $\mathcal{N}$ of bdry$\mathcal{K}$.
\item[(A9)]\quad For all $(x,p) \in \mathbb{R}^n \times \left( \mathbb{R}^n\setminus \{0\} \right)$, there exists a unique $u^* \in U$ such that
\begin{equation}
-f(x,u^*).p - L(x,u^*) = \mathcal{H}(x,p),
\end{equation}
and the function $u^*: (x,p) \mapsto u^*(x,p)$ is of class $C^k$ in $\mathbb{R}^n \times \left( \mathbb{R}^n\setminus \{0\} \right)$.
\end{itemize}

For our analysis, in assumption (A9), we only need that $u^*$ is of class $C^k$ in an open neighborhood of the set $\{(x,p)\in \mathbb{R}^n \times \left(\mathbb{R}^n\setminus \{0\}\right): \mathcal{H}(x,p) = 0\}$. For examples satisfying this condition, one can find in \cite{CPS00,P02}. It follows from (A8) and (A9) that the Hamiltonian satisfies
\begin{itemize}
\item[(H3)]\quad $\mathcal{H} \in C^k\left(\mathbb{R}^n \times \left( \mathbb{R}^n\setminus \{0\} \right)\right)$.
\end{itemize} 
We next introduce the definition of conjugate times which is related to the Jacobians of solutions of the backward Hamiltonian system considered in \cite{P02}. Given $z\in \mathrm{bdry}\mathcal{K}$, we denote by $(y(z,\cdot),q(z,\cdot))$  the solution the \textit{backward Hamiltonian system}
  \begin{equation}
\label{BHS}
\left\{
  \begin{array}{lcl}
    \dot{y}(t) &=& \mathcal{H}_p(y(t),q(t))\\
  - \dot{q}(t) &=& \mathcal{H}_x(y(t),q(t))
  \end{array}
  \right.
\end{equation}
with the initinal conditions
\begin{equation}
\label{IC}
\left\{
  \begin{array}{lcl}
    y(0) &=& z\\
   q(0) &=& \varphi(z),
  \end{array}
  \right.
\end{equation}
where $\varphi(z) = \nabla \psi(z) +\mu (z) \nabla b_\mathcal{K}(z)$ with $\mu(\cdot)$ satisfying 
$$\mathcal{H}(z,\nabla \psi(z) + \mu(z)\nabla b_{\mathcal{K}}(z)) = 0.$$ 
Note that $\mathcal{H}(z,\varphi(z)) = 0$ and that under our assumptions the function $\varphi: \mathrm{bdry}\mathcal{K} \to \mathbb{R}^n$ is of class $C^k$ (see, e.g., \cite{P02}).

As shown in \cite{P02}, the solution $(y(z,\cdot), q(z,\cdot))$ of (\ref{BHS}) - (\ref{IC}) is defined for all $t\in [0,+\infty)$. Moreover, $y(\cdot), q(\cdot)$ are of class $C^k$ on $\mathrm{bdry}\mathcal{K} \times [0,+\infty)$.

Now let $Y$ and $Q$ denote respectively the Jacobians of  $y(\cdot)$ and $q(\cdot)$ with respect to the pair $(z,t)$ in $\mathrm{bdry}\mathcal{K} \times [0,\infty)$ where $(y(z,\cdot),q(z,\cdot))$ solves (\ref{BHS}) - (\ref{IC}). Then $(Y,Q)$ is the solution of the system
 \begin{equation}
\label{VS}
\left\{
  \begin{array}{lcl}
    \dot{Y} &=& \mathcal{H}_{xp}(y(z,t),q(z,t))Y + \mathcal{H}_{pp}(y(z,t),q(z,t))Q\\
   -\dot{Q} &=& \mathcal{H}_{xx}(y(z,t),q(z,t))Y + \mathcal{H}_{xp}(y(z,t),q(z,t))Q,
  \end{array}
  \right.
\end{equation}
with the initial conditions
 \begin{equation}
\label{VSIC}
\left\{
  \begin{array}{lcl}
    Y(z,0) &=& A(z)\\
   Q(z,0) &=&  B(z),
  \end{array}
  \right.
\end{equation}
where $A(z), B(z)$ are square matrices depending smoothly on $z$ which we can compute.

As explained in \cite{P02}, the Jacobian $Y$ and $Q$ are understood in the following sense. Fixed $z_0 \in \mathrm{bdry}\mathcal{K}$ and $t_0 >0$. Since $\mathrm{bdry}\mathcal{K}$ is an $(n-1)$-dimensional manifold of class $C^{k+1}$, there exist an open neighborhood $I$  and a parameterized function $\xi: z\in I \mapsto \eta \in \xi(I) \subset \mathbb{R}^{n-1}$ of class $C^{k+1}$ with the inverse $\phi$ of class $C^{k+1}$, where $\psi(I)$ is an open neighborhood of $\eta_0 = \xi (z_0)$.  Then $Y(z_0,t_0)$ and $Q(z_0,t_0)$  denote the Jacobians of $Y(\phi(\cdot),\cdot)$ and $Q(\phi(\cdot),\cdot)$ with respect to the coordinates $\eta \in  \mathbb{R}^{n-1}$ and the time $t$ at the point $(\eta_0,t_0)$ i.e., $Dy(\phi(\eta_0),t_0)$ and $Dq(\phi(\eta_0),t_0)$. In this case, 
$$A(z_0) = Y(z_0,0) = \left(\mathcal{H}_p(z_0,\varphi(z_0)), \frac{\partial\phi}{\partial \eta}(\eta_0)\right),$$
and one can compute that $\det A(z_0) = \alpha  \mathcal{H}_p(z_0,\varphi(z_0)).\varphi(z_0)$ for some real constant $\alpha \ne 0$. Therefore, $\det A(z_0) = \det Y(z_0,0) \ne 0$ (see proof of Lemma 4.2 in \cite{P02}). Then
$$\mathrm{rank}\begin{pmatrix}
Y(z_0,0)\\ Q(z_0,0)
\end{pmatrix} = n
$$
and by properties of linear systems, we have
$$
\mathrm{rank}\begin{pmatrix}
Y(z_0,t)\\ Q(z_0,t)
\end{pmatrix} = \mathrm{rank}\begin{pmatrix}
Y(z_0,0)\\ Q(z_0,0)
\end{pmatrix} = n,\quad \forall t\in [0,+\infty).
$$
Note that this definition of the Jacobian depends on the parameterized function. For our purpose, however, this does not matter because we only focus on the ranks of the matrices $Y$ and $Q$ which are independent of the choice of the parameterized functions.
\begin{definition}
For $z\in\mathrm{bdry}\mathcal{K}$, the time
$$t_c(z): = \inf\{t\in [0,+\infty): \det Y(z,s) \ne 0,\,\forall s\in [0,t]\}$$
is said to be conjugate for $z$ if and only if $$t_c(z) <+\infty\quad \text{and} \quad \det Y(z,t_c(z)) = 0.$$
\end{definition}
 
 Fixed $z_0\in \mathrm{bdry}\mathcal{K}$ and $T_0>0$. Since $Y(z_0,t)$ is invetible for $t$ sufficiently small, if there exists conjugate time $t_c$ for $z_0$ then $t_c >0$. On the other hand, if there is no conjugate time for $z_0$ in $[0,T_0]$, then by the continuity and the fact that $\det Y(t,z_0) \ne 0$ for all $t\in [0,T_0]$, there exist $\varepsilon >0, \sigma >0$ such that there is no conjugate time for any $z\in B(z_0,\sigma)\cap \mathrm{bdry}\mathcal{K}$ in $[0,T]$ with $T< T_0+\varepsilon$. In this case $y(\cdot)$ is an one-to-one correspondence in a neighborhood of $(z_0,T_0)$. Using this fact, one can prove the following theorem.
 
 \begin{theorem}\label{M1}
Assume (A0) - (A6) and (A8) - (A9). Let $\bar{x}\in \mathcal{R}\setminus \mathcal{K}$ be such that $V$ is differentiable at $\bar{x}$. Let $x(\cdot)$ be the optimal trajectory starting at $\bar{x}$ and $\tau$ be the exit time of $x(\cdot)$. Set $\bar{z}= x(\tau)\in \mathrm{bdry}\mathcal{K}$. If there is no conjugate time for $\bar{z}$ in $[0,\tau]$ then $V$ is of class $C^k$ on an open neighborhood of $\bar{x}$.
 \end{theorem}
 
 \begin{remark}
 Observe that if $V$ is differentiable at a point $\bar{x}\in \mathcal{R}\setminus \mathcal{K}$ then V is differentiable along the optimal trajectory starting at $\bar{x}$ except the final point. Then in Theorem \ref{M1} we can conclude that $V$ is of class $C^k$ on an open neighborhood of $x(s)$ for all $s\in [0,\tau)$. 
 \end{remark}
 
Following the idea used in \cite{CFS14} where the authors study the regularity of the value function for a Mayer optimal control problem, by using Theorem \ref{M1}, we can prove the following theorem.
\begin{theorem}\label{M2}
AsAssume (A0) - (A6) and (A8) - (A9).  Let $\bar{x} \in \mathcal{R}\setminus \mathcal{K}$. If $\partial^PV(\bar{x}) \ne \emptyset$ then $V$ is of class $C^k$ on an open neighborhood of $\bar{x}$. 
\end{theorem}

\begin{remark} In the case the running cost does not depend on $u$- variable, i.e., $L = L(x)$, the results of this section still hold true if we repace assumption (A9) by (H3), a weaker assumption.
\end{remark}

\textbf{Acknowledgments}. This chapter was started writing when the author was a SADCO PhD fellow, position ESR4 at University of Padua which was supported by the European Union under the 7th Framework Programme “FP7-PEOPLE-2010-IT”, Grant agreement number 264735-SADCO. This work was also supported by funds allocated to the implementation of the international co-funded project in the years 2014-2018, 3038/7.PR/2014/2, and by the EU grant PCOFUND-GA-2012- 600415.





\end{document}